\documentstyle[aps,preprint,pra,tighten]{revtex}
\begin{document}
\draft
\title{An Adiabatic Theorem without a Gap Condition: Two level system
coupled to quantized radiation field}
\author{J.~E.~Avron and A.~Elgart}
\address{Department of Physics, Technion, 32000 Haifa, Israel}
\maketitle
\begin{abstract} We prove an adiabatic theorem for the ground state of the
Dicke model in a slowly rotating magnetic field and show that for weak
electron-photon coupling, the adiabatic time scale is close to the time
scale of the corresponding two level system--without the quantized
radiation field. There is a correction to this time scale which is the Lamb
shift of the model. The photon field  affects the  rate of approach to the
adiabatic limit through a logarithmic correction originating from an
infrared singularity characteristic of QED.
\end{abstract}
\section{Introduction} In this work we investigate  the relation between
adiabatic theorems for models that, like QED,  allow for the creation and
annihilation of photons, and the  corresponding quantum mechanical models
where the electron is  decoupled from the photon field. We study this
problem in the context of a specific and  essentially soluble model: The
Dicke model
\cite{d}. The corresponding quantum mechanical model is a two level system,
such as a spin in an adiabatically rotating magnetic field,  which is a
basic paradigm of adiabatic theory \cite{berry}.

In the usual quantum adiabatic theorem \cite{bf,kato} the gap between
eigenvalues plays an important role:  It fixes the adiabatic time scale
and  determines the rate at which the adiabatic limit is approached.
 There is no such gap in  the corresponding QED models so the nature of the
adiabatic theorem in the two cases has qualitatively different features.
For example, there is no gap in the spin-boson and Dicke models (for weak
coupling) both of which describe a two level system in a radiation field.

The first problem we address is whether there is an adiabatic theorem for
the ground state in a radiation field. Assuming a positive answer, the
second question is, what property of the QED model, plays the role of the
gap in the adiabatic theorem. Another way of phrasing this question is how
does the adiabatic time scale  of the two level system compare to that of
the QED model? Are the two  close in the limit of small fine structure
constant, $\alpha$, and if so, how close? The third question compares the
rate of approach to the adiabatic limit in the two models.

Consider a two level system, such as a spin or a twofold Zeeman split
atomic  level, in an external magnetic  field pointing in the z direction.
When  radiation effects are neglected, the  corresponding Hamiltonian is
\begin{equation} H= m\sigma_z,\quad m=\mu B .
\end{equation}
 The corresponding Dicke model is
\begin{equation} H_D=H\otimes {\bf 1} + \alpha^{-1}\, {\bf 1}\otimes E +
\sqrt{\alpha}\,\sigma_+\otimes a^{\dagger}(f) +
\sqrt{\alpha}\,\sigma_-\otimes a(f),
\end{equation} where
\begin{equation} E=\int{|k|\, a^{\dagger}(k)a(k)d^dk},\end{equation}
and\footnote {$\{\cdot,\cdot\}$ stands for anticommutator}
\begin{equation} f(k)=\,\sqrt{\frac{2\pi}{|k|}}
\,\langle\psi_1\vert\,\{e^{-i{k}\cdot {x}}, p\}\,\vert\psi_2\rangle,
\label{fk}
\end{equation} with $\psi_j(x)$ the atomic wave functions of the  two
level system.  Polarization indices are omitted  since the helicity of the
photon does not play an interesting role in the questions we study. We use
atomic units where
$e=\hbar=1$ so $\alpha=1/c=1/137$ is small.  $\mu$, the magnetic moment,
is also of order $\alpha$ in these units.

Following M. Berry's  \cite{berry} let us consider the case where the
magnetic field changes its direction adiabatically and has fixed
magnitude. The time dependent  Hamiltonian for the two level system is
\begin{equation} H(s)=\mu B(s)\cdot\sigma = U(s)HU^*(s),
\end{equation} with $U(s)\in SU(2)$ the appropriate rotation. The
corresponding adiabatic Dicke model has the time dependent Hamiltonian:
\begin{equation} H_D(s)=\left(U(s)\otimes {\bf 1}\right)\, H_D\,
\left(U^*(s)\otimes {\bf 1}\right).
\end{equation} Our aim is to compare the adiabatic evolution of the ground
state of $H(s)$ with that of $H_D(s)$.

Adiabatic theorems for quantum systems coupled to a field have been
studied in
\cite{nt,ds}. In \cite{nt} Narnhofer and Thirring give characterization of
extremal KMS states by  adiabatic invariance. When applicable, this result
says, in particular, that the ground state is adiabatic invariant. The
 characterization depends on a condition of asymptotic abelianess which
does not hold for the models we consider. In \cite{ds} Davis and Spohn
give a  derivation  of  linear response theory for a system coupled to a
bath in the adiabatic limit. The notion of adiabaticity in this work  is
such that the coupling between the field and the quantum system vanishes
in the adiabatic limit.  This is not a standard notion of adiabaticity.

Let us now describe our results. First, we show that there is an adiabatic
theorem for the ground state of the Dicke model, even though the model has
no spectral gap to protect the ground state. Second, we show that the the
distance to a nearby resonance in the Dicke model plays the role of a gap.
Third, we show that the adiabatic time scale for the Dicke model and the
two level system agree in the limit of small $\alpha$. The (inverse) of
two time scales differ by the Lamb shift of the Dicke model. And finally,
we show that the approach to the adiabatic limit in the two models is
different: While in the two level system the approach to the adiabatic
limit is with an error $O\left(\frac{1}{\tau}\right)$, the approach to the
limit in the Dicke model is with an error of
$O\left(\frac{\sqrt{\log \tau}}{\tau}\right)$. The logarithm comes from an
infrared divergence characteristic of QED.

Although the results we derive here are for a rather special model we
suggest that something similar  happens also for more realistic models.
The success of the quantum  in numerous applications that depend on a
correct prediction of the adiabatic time scale is evidence that at least
the time scale aspect of our results may well carry over to more realistic
models. It would be interesting to know if this is indeed the case for the
Spin-Boson model \cite{hs,dg,ms}. The spin-boson model is a more realistic
QED version of a two level system, which, unlike the Dicke model, is not
explicitly soluble. However, as much progress in the spectral analysis of
the spin-boson problem has been recently made, the problem we pose here
may be a reasonable challenge.

\section{The Adiabatic Theorem and A Commutator Equation}

In this section we  explain what we mean by ``adiabatic theorem'', and
give a condition for an adiabatic theorem to hold. This condition is that
the commutator equation, Eq.~(\ref{commutators}) below, has solutions $X$,
$Y$ which are bounded operators\footnote{for $X$ we also need that its
derivative is bounded}. We also introduce notation, terminology, and
collect known facts that we need. To simplify the presentation, we shall
stay away from making optimal assertions.

We consider Hamiltonians that are bounded from below, and choose the
origin of the energy axis so that  the spectrum begins at zero. Let
$H(s)\ge 0$ be a family of such self-adjoint Hamiltonians.  The unitary
evolution generated by the Hamiltonian, $U_\tau(s)$, is the solution of
the initial value problem:
\begin{equation} i\,\dot U_\tau (s) = \tau H(s) U_\tau(s),\quad
U_\tau(0)=1,
\quad s\in [0,1].\label{schrodinger}
\end{equation}
$\tau$ is the adiabatic time scale, and we are concerned with the limit of
large $\tau$. The physical time is $t=\tau s\in[0,\tau]$.  Since
$\tau$ is large $H(s)=H(t/\tau)$  varies adiabatically. We assume that all
operators are defined on some fixed dense domain in the Hilbert space.

The (instantaneous) ground state is in the range of the kernel of $H(s)$
and we assume that the kernel is smooth and  one-dimensional. Let
$P(s)\neq 0$ be the projection on the kernel of $H(s)$, i.e. $H(s)\,
P(s)=0$, $dim\,P= Tr\,P=1$. By smoothness we mean that $\dot P(s)$ a
bounded operator.

The adiabatic theorems we consider are  concerned with  the large time
behavior of the evolution of the ground state where $t=O(\tau)$ or,
equivalently,$s=O(1)$. The smoothness of the kernel implies that there is
a natural candidate for an
 adiabatic theorem for the ground state, which is independent of whether
$H(s)$ does or does not have a gap in it spectrum.  Namely, that if
$\psi(0)\in Range \, P(0)$ at time $s=0$, then it evolves in time so that,
$\psi_\tau(s)=U_\tau (s)\ \psi(0)$ lies in  $Range\, P(s)$ at time $s$  in
the adiabatic limit, $\tau\to\infty$.

To formulate the adiabatic theorem with error estimates we need to get
hold of  {\em adiabatic phases} \cite{berry}. To do that we introduce the
adiabatic evolution of Kato \cite{kato}: Let $U_A(s)$ be the solution of
the evolution equation
\begin{equation} \dot U_A(s) = [\dot P(s),P(s)]\, U_A(s),\quad U_A(0)=1,
\quad s\in [0,1].\label{kato}
\end{equation} It is known that
\begin{equation}  U_A(s)\, P(0) =P(s)\, U_A(s).
\end{equation} That is $U_A(s)$ maps $Range\ P(0)$ onto $Range\ P(s)$. We
can now formulate the basic adiabatic theorem :
\\\\
{\bf {Theorem\  II.1}}
 Let $H(s)P(s)=0$ for all $0\le s\le 1$, with $P$
differentiable  projection on the ground state, with $\Vert
\dot P(s)\Vert\le D$. Suppose that the commutator equation
\begin{equation} [\dot P(s), P(s)] = [H(s),X(s)] +Y(s),\label{commutators}
\end{equation} has operator valued solutions, $X(s)$ and $Y(s)$ so that for
$\varepsilon\searrow 0$
\begin{equation}
\Vert X(s)\Vert\ + \Vert \dot X(s)\Vert\ \le\ C\left\{ \begin{array}{c}
\,\varepsilon^{-\nu}\\ \vert\log \varepsilon\vert
\end{array} \right. ,\quad \Vert Y(s)\Vert \le \hat C\,\varepsilon^{\mu},
\end{equation} with $\mu,\nu\ge 0$. Then
\begin{equation}
\Vert (U_\tau (s)-U_A(s))P(0)\Vert \le \tilde C \left\{
\begin{array}{c}\tau^{-\frac{\mu}{\nu+\mu}}
\\ \frac{\log \tau}{\tau}\end{array}\right. ,\quad s\in[0,1].
\end{equation}

 Remarks: 1. In the case that there is a gap in the
spectrum, one can always find  $X(s)$ bounded so $\nu=0$,  and $Y=0$,  see
\cite{asy}. $X$, and therefor also $\tilde C$, is of the order of
(gap)$^{-1}$.  This gives  error of $1/\tau$, and generalizes the
adiabatic theorem of Born and Fock and Kato  for discrete spectra, to more
complicated spectra provided there is a gap.

2. The theorem says that the physical evolution clings to the instantaneous
spectral subspace. In particular, if $P$ is one dimensional, it says that
the physical evolution of the ground state remains close to the
instantaneous ground state.

3. Here, and throughout, we are concerned only with the adiabatic theorem
to
 lowest order. If $s$ is chosen outside the support of $\dot P$ then much
stronger results can be obtained. See e.g. \cite{ks}.

4. The adiabatic time scale $\tau_0$ set by this theorem is $\tau_0=
O((2+D)C)$.

Proof:  Let $W(s) =  U_A^\dagger(s) U_\tau(s))$, with $W(0)=1$. From the
equation of motion, and the commutator equation, Eq.~(\ref{commutators}),
\begin{eqnarray}
 P(0)\,\dot W(s)\,&=& -P(0) U^\dagger_A(s)\Big(i\,\tau\,H(s)+[\dot
P(s),P(s)]\Big) U_\tau(s)\nonumber
\\ &=&-U^\dagger_A(s)\,P(s)\Big(i\,\tau\,H(s)+[\dot
P(s),P(s)]\Big)U_\tau(s)
\nonumber
\\ &=&\,-U^\dagger_A(s)\,P(s)\,[\dot P(s),P(s)]\, U_\tau(s)
\nonumber
\\ &=&\,-U^\dagger_A(s)\,P(s)\,\Big([H(s),X(s)]+Y(s)\Big)\,U_\tau(s)
\nonumber
\\ &=&\,-U^\dagger_A(s)\,P(s)\,\Big(-X(s)\,H(s)+Y(s)\Big)\,U_\tau(s)
\nonumber
\\ &=&\frac{i}{\tau}\,P(0)\, U^\dagger_A(s)\,X(s)\,
\dot U_\tau(s)-P(0)\,U^\dagger_A(s)\,Y(s)\, U_\tau(s).
\end{eqnarray} To get rid of derivatives of $U_\tau$, which are large by
the equation of motion, we rewrite the first term on the rhs (up to the
$P(0)$ on the right) as :
\begin{eqnarray}
 U^\dagger_A(s)\,X(s)\,\dot U_\tau(s) &=&\dot{\left( U^\dagger_A(s)\,X(s)
U_\tau(s)\right)}- U^\dagger_A(s)\,\dot X(s)\, U_\tau(s)-
\dot U^\dagger_A(s)\,X(s)\,U_\tau(s) \\&=&
\dot{\left( U^\dagger_A(s)\,X(s)\, U_\tau(s)\right)} - U^\dagger_A(s)\,\dot
X(s)\, U_\tau(s)+ U^\dagger_A(s)\,[\dot P(s),P(s)]\,X(s)\,
U_\tau(s).\nonumber
\end{eqnarray} From this it follows, by integrating,  that for $s\in[0,1]$
\begin{eqnarray}
\Vert (U_\tau(s)-U_A(s))\,P(0)\Vert &=& \Vert
P(0)(U^\dagger_\tau(s)-U^\dagger_A(s))\,\Vert\nonumber \\
\Vert P(0)(1-W(s))\,\Vert&\le& \hat C\varepsilon^\mu +\frac{(2+D)\,
C}{\tau}\left\{
\begin{array}{c}{\varepsilon^{-\nu}}
\\ \vert\log\varepsilon\vert.\end{array}\right.
\end{eqnarray} Choosing $\varepsilon= \tau^{-\frac{1}{\mu+\nu}}$ gives
\begin{equation}
\Vert (U_\tau(s)-U_A(s))\,P(0)\Vert  \le \tilde C \left\{
\begin{array}{c}\tau^{-\frac{\mu}{\nu+\mu}}
\\ \frac{\log \tau}{\tau}\end{array}\right. .
\end{equation}
 This concludes the proof of the theorem.\hfill$\hfill\vbox{\hrule height
0.6pt
  \hbox{\vrule width 0.6pt height 1.8ex \kern 1.8ex\vrule width 0.6pt}
        \hrule height 0.6pt}$

It is convenient to rewrite this solvability condition in a way that one
needs to solve for a fixed $X$ and $Y$ rather than functions $X(s)$ and
$Y(s)$. This is accomplished by
\\\\
{\bf {Corollary\  II.1.1}}
 Let $P(s)$ be the family
\begin{equation} P(s)=V(s)\, P\,V^\dagger(s),\quad V(s)=\exp
(i\,s\,\sigma).
\end{equation} It is  enough to solve for the commutator equation
\begin{equation}
 i\, K=[H, X] +Y, \quad K= \{\sigma,P\} -2 P\sigma P,\label{commute}
\end{equation} for fixed $X$ and $Y$ so that for $\varepsilon\searrow 0$
\begin{equation}
\Vert X\Vert\ \le\ C\left\{ \begin{array}{c}
\,\varepsilon^{-\nu}\\ \vert\log \varepsilon\vert
\end{array} \right. ,\quad \Vert Y\Vert
\le \hat C\,\varepsilon^{\mu},
\end{equation} with $\mu,\nu\ge 0$, and   $\Vert \dot P(s)\Vert\le D$.
$X(s)$ and
$Y(s)$ are then determined by the obvious unitary conjugation.

 Proof:
 Since $P(s)=V(s)\, P\,V^\dagger(s)$, we have
\begin{equation}
\dot P(s)=i\,V(s)\, [\sigma,P]\,V^\dagger(s),
\end{equation} and
\begin{eqnarray} [\dot P(s),P(s)]&=&i\,V(s)\,
\Big[[\sigma,P],P\Big]\,V^\dagger(s)\nonumber \\ &=&i\,V(s)\,
\Big(\{\sigma,P\}-2P\,\sigma\,P\Big)\,V^\dagger(s).
\end{eqnarray}
\hfill$\hfill\vbox{\hrule height 0.6pt\hbox{\vrule width 0.6pt height 1.8ex
\kern 1.8ex \vrule width 0.6pt}\hrule height 0.6pt}$

\section{ An Adiabatic Theorem for a Threshold State: The Friedrichs Model}
As a warmup, and as a preparation for the analysis of the Dicke model, let
us prove an adiabatic theorem for the Friedrichs model which has a bound
state at the threshold of the continuum.

There is an inherent difficulty in the situation of a bound state at
threshold in general,  and in the Friedrichs model
\cite{friedrichs,friedrichs1,howland} in particular, namely, that a bound
state at threshold is not a stable situation. Under a small deformation of
the Hamiltonian, the ground state will, generically, split away from the
absolutely continuous spectrum and a gap  develops. Since our aim is to
study families related by a unitary, this problem does not appear. That
is, we consider the family
$H_F(s)=V(s)\,H_F\,V^\dagger(s)$ where $H_F$ has a bound state at
threshold and
$V(s)$ is a smooth family of unitaries.

\subsection{The Friedrichs Model}

We shall consider a family of Hamiltonians, closely related to the standard
Friedrichs model \cite{friedrichs}, parameterized by the scaled time $s$,
a real number $d>0$ that plays the role of dimension, and a function $f$
that describes the deformation of the family. Since we are only interested
in the low energy behavior of the family we shall introduce an
``ultraviolet cutoff'' to avoid inessential difficulties.

The Hilbert space of the  Friedrichs model (with an ultraviolet cutoff) is
${\cal H}=\kern.1em{\raise.47ex\hbox{  $\scriptscriptstyle
|$}}\kern-.40em{\rm C}\oplus L^2([0,1],k^{d-1}\,dk)$. A vector
$\psi\in{\cal H}$ is normalized by
\begin{equation}
\psi=\left(\begin{array}{c}\beta\\ f(k)
\end{array}\right)\, \quad \Vert \psi\Vert^2= |\beta\vert^2+\int_0^1 \vert
f(k)\vert^2 k^{d-1}\,dk,\quad
\beta\in\kern.1em{\raise.47ex\hbox{  $\scriptscriptstyle
|$}}\kern-.40em{\rm C}.
\end{equation} We choose a special, and trivial, case of a diagonal
Hamiltonian whose action on a vector $\psi$ is as follows:
\begin{equation} H_F\,\psi=\left(\begin{array}{ll} 0&0\\ 0&k
\end{array}\right)\,\left(\begin{array}{c}\beta\\ f(k)
\end{array}\right)=\left(\begin{array}{c} 0\\ k\,f(k)
\end{array}\right).
\end{equation}
 $H$ has a ground state at zero energy with projection
\begin{equation} P=\left(\begin{array}{ll} 1&0\\ 0&0
\end{array}\right).
\end{equation} The rest of the spectrum is the unit interval $[0,1]$, and
is absolutely continuous. The density of states in this model is
proportional to
$E^{d-1}$.

We construct the family $H(s)$ by conjugating $H$ with a family of
unitaries:
\begin{equation} V_f(s)=\exp is\sigma(f), \quad
\sigma(f)=\left(\begin{array}{ll} 0&\langle f\vert \\ \vert f \rangle&0
\end{array}\right),
\end{equation} where $f$ is a  vector in $L^2([0,1],k^{d-1}\, dk)$.
\\\\
{\bf {Theorem\  III.1}}
 Let $H_F(s;d,f)$ be the family of Friedrichs models with
a ground state at threshold for all s
\begin{equation} H_F(s;d,f)= V_f(s)\, H_F\,V^\dagger_f(s).
\end{equation}
 Suppose that
\begin{equation} g(k)=i\, k^{-1}\, f(k)\in L^2([0,1],k^{d-1} dk), \quad
V_f(s)=\exp\, i\,s\,\sigma(f)\, ,
\end{equation} then the quantum evolution of the ground state of
$H_F(s;d,f)$ is adiabatic and its deviation from the instantaneous ground
state is, at most,
$O(1/\tau)$.

 Remarks: 1. Note that if the conditions in the theorem
hold in dimension $d_0$, then they hold in all dimensions $d\ge d_0$. The
physical interpretation of that is that the density of states at low
energies decreases with $d$. So, even though there is spectrum near zero,
there is only very little of it.

2. If $g$ is not in $L^2$ there may still be an adiabatic theorem with
slower falloff in $\tau$ by accommodating $Y\neq 0$. An example will be
discussed in the next section.

3. The Friedrichs model is vanilla: $H_F$ has no interesting energy scale
to fix the adiabatic time scale. The scale is set by the perturbation
alone:
$\tau_0= O((1+\Vert f\Vert^2)\,\Vert g\Vert )$. This is quite unlike the
case in the usual adiabatic theorem and unlike what we shall show for the
Dicke model.

Proof: In this case $K$ of Corollary 2.1 is  $K=\sigma(f)$. With $g\in
L^2$,
$\sigma(g)$ is a bounded (in fact, finite rank) operator and an easy
calculation gives
\begin{equation} [H_F, \sigma(g)]=\left(\begin{array}{lr}0&\langle -kg|\\
|kg\rangle&0\end{array}\right)=i\sigma(f).
\end{equation} Hence  \begin{equation} X= \sigma(g) ,\quad Y=0\, ,
\end{equation} solve the commutator equation, Eq.~(\ref{commute}), with a
bounded
$X(s)$ and $Y(s)=0$.\hfill$\hfill\vbox{\hrule height 0.6pt\hbox{\vrule
width 0.6pt height 1.8ex \kern 1.8ex \vrule width 0.6pt} \hrule height
0.6pt}$


\section{Adiabatic Theorem for the Dicke Model}

In this section we describe an adiabatic theorem for the Dicke model
\cite{d} that says that the an adiabatic rotation of a two level system
evolves the ground state so that it adheres to the instantaneous ground
state and the time scale, at least in three dimensions, is essentially the
time scale fixed by Quantum Mechanics without photons. The rate of
approach to the adiabatic limit is different from that of a two level
system and has a logarithmic correction in three dimensions. This section
also collects known facts about the Dicke model that we need.

\subsection{The Dicke Model}

The Spin-Boson Hamiltonian in the canonical QED version of a two level
system
\cite{hs,dg,ms}. The Dicke model is a simplified version of the Spin-Boson
Hamiltonian in the rotating wave approximation. The rotating wave
approximation,  can indeed be motivated in the single-mode Dicke model.
In  the multi-mode case we consider the rotating wave approximation is  a
name that describes which terms in the Spin Boson Hamiltonian are kept and
which are not.

The model describes a two level system coupled to a massless boson field
in $d$ dimensions.  The Hamiltonian is:
\begin{equation} H_D(m,d,f,\alpha)=m\,(1 -P)\otimes {\bf 1} + \alpha^{-1}\,
{\bf 1}\otimes E +
\sqrt{\alpha}\,\sigma_+\otimes a^{\dagger}(f) +
\sqrt{\alpha}\,\sigma_-\otimes a(f),\label{eq1}
\end{equation}acting on the Hilbert space $\kern.1em{\raise.47ex\hbox{  $
\scriptscriptstyle|$}}\kern-.40em{\rm C}^2\otimes {\cal F}$ with ${\cal
F}$ being the symmetric Fock space over $L^2(R^d,d^dk)$. Here
\begin{equation} P=\left(\begin{array}{ll} 1&0\\
0&0\end{array}\right),\quad
\sigma_+=\left(\begin{array}{ll} 0&1\\ 0&0\end{array}\right),
\quad
\sigma_-=\left(\begin{array}{ll} 0&0\\ 1&0\end{array}\right)\quad
E=\int{|k|\, a^{\dagger}(k)a(k)d^dk}.
\end{equation}
$m>0$ is the gap in the quantum Hamiltonian (without photons). $a(f)$ and
$a^\dagger(f)$ are the usual creation and annihilation operators on ${\cal
F}$ obeying the canonical commutation  relations
\begin{equation} [a(f),a^{\dagger}(g)]=\langle f\vert g\rangle.
\end{equation} We denote by $|0\rangle$ the field vacuum and by $\Omega$
the projection on the vacuum.

It may be worthwhile to explain where the various powers of $\alpha$ in
$H$ come from. For the radiation field  the $\alpha^{-1}$ comes from
$\hbar\omega =\hbar c |k|$  which explains why the field energy comes with
a large coupling constant. The $\sqrt{\alpha}$  has one inverse power of
$c$  from  minimal coupling, $\frac{e}{2mc} (p\cdot A+A\cdot p)$. Half a
power of $\alpha$ comes from the standard formula for  the vector potential
\begin{equation} {A}({x}):=
\int d^3k\,\sqrt{\frac{2\pi c}{|k|}}\,\Big(e^{-i{k}\cdot
{x}}\,a^{\dagger}({k}) +e^{i{k}\cdot{x}} a({k})\Big).
\end{equation} Compare e.g. \cite{dicke}.

With reasonable atomic eigenfunctions, $f(k)$, Eq.~(\ref{fk}) has fast
decay at infinity and the model is ultraviolet regular. In the infrared
limit $f(k)$ behaves like
\begin{equation} f(k)\to-i\,\sqrt{\frac{2\pi}{|k|}}
\int\,\Big(\psi_1^*(x)\,(\nabla\psi_2)(x)\,
-(\nabla\psi_1)^*(x)\,\psi_2(x)\Big)\,d^dx.
\end{equation} In particular we see that for small $k$
 \begin{equation} f(k)=K\,\sqrt{\frac{1}{ |k|}}.\label{K}
\end{equation} The square root singularity is a characteristic infrared
divergence of QED, and it has consequences for the adiabatic theorem as we
shall see. Note that with $f$ having a square root singularity the model
makes sense (as an operator) provide $d>1$, for otherwise $a^\dagger(f)$ is
ill defined since $f$ is not in $L^2$.

An important parameter in the model is
\begin{equation} {\cal E}= \left\langle f \left\vert
\frac{1}{|k|}\right\vert f\right\rangle.
\end{equation} Bearing in mind the square root singularity  of $f$ we see
that
\begin{equation} {\cal E}\sim \int \frac{d^dk}{|k|^2},
\end{equation} is  finite for all $d>2$.


\subsection {Spectral Properties} What makes the Dicke model simple is
that it has a constant of motion \cite{hs}. If we let $N=\int a^\dagger(k)
\, a(k) d^dk$ be the photon number operator, then ${\cal N}$ commutes with
$H_D$ where
 \begin{equation} {\cal N}=\left(\begin{array}{ll} N&0 \\
0&N+1\end{array}\right)={\bf 1}\otimes N + P\otimes {\bf 1}.
\end{equation} The spectrum of ${\cal N}$ is the non-negative integers. The
spectral properties of $H_D(m,d,f,\alpha)$ can be studied by restricting to
subspaces of ${\cal N}$.

\paragraph{${\cal N}=0$}: The kernel of ${\cal N}$ is one dimensional and
is associate with the projection
\begin{equation} P=\left(\begin{array}{ll}\Omega&0 \\
0&0\end{array}\right).
\end{equation}
$\Omega$ is the projection on the field vacuum. It is easy to see that
$P\,H_D(m,d,f,\alpha)P=0$, so the model always has a state at zero energy.
This state may or may not be the ground state. It is the ground state if
$\alpha^2 {\cal E} <m$ \cite{hs}.

\paragraph{${\cal N}=1$}: The space is  basically ${\cal H}$ of the
Friedrichs model. The correspondence of vectors in the two spaces is
\begin{equation} \left(\begin{array}{c}a^\dagger(g)\\ \beta
\end{array}\right)\,|0\rangle\leftrightarrow  \left(\begin{array}{c}g\\
\beta
\end{array}\right).
\end{equation}
 The Hamiltonian action in the Friedrichs model language is:
\begin{equation}
H_D(m,d,f,\alpha)\leftrightarrow\left(\begin{array}{ll}\frac{|k|}{\alpha}&|
\sqrt{\alpha}\,f\rangle\\ \langle \sqrt{\alpha}\, f|&m\end{array}\right)\,.
\end{equation} It is a standard fact about the Friedrichs model
\cite{friedrichs,friedrichs1} that provided
\begin{equation}
\alpha^2 {\cal E} <m,\label{standard}
\end{equation} the model has no bound state, and the spectrum is
$[0,\infty)$ and is absolutely continuous. Since $f$ has square root
singularity at the origin, (and has fast decay at infinity), this
condition holds for $d\ge 3$ if
$\alpha $ (or $f$)  is not too large. In three dimensions, provided that
the level spacing $m>>\alpha^2 $ in atomic units, (about $10^{-3}$ eV), the
inequality  holds. In two dimensions the left hand side is log divergent,
and the spectrum in the ${\cal  N}=1$ sector has a bound state at negative
energy. This state lies below the bound state of the ${\cal  N}=0$ sector.
We do not  consider this situation and henceforth stick to $d\ge 3$.

\paragraph{${\cal N}\ge 2$}: It is known \cite{hs} that the bottom of the
spectrum in all these sectors is at zero if (\ref{standard}) holds.

\subsection{Adiabatic Rotations}

Suppose, that the two level system of the Dicke model describes e.g. two
Zeeman split energy levels of an atom in constant external magnetic field
$B$ pointing in the z direction.  Rotations about the z axis do not change
the orientation of the magnetic field, and commute with $\cal N$ and are
uninteresting. Rotations about the x axis change the orientation of the
magnetic field and are implemented by
\begin{equation} V(s) =\exp\, (i\,s\, \sigma_x)\otimes {\bf 1}.\label{v}
\end{equation} Such rotations do not commute with ${\cal N}$. Indeed,
\begin{equation}[{\cal N},\sigma]=\left(\begin{array}{rl}0&{ 1}\\
-{1}&0\end{array}\right)\otimes {\bf 1}=J\otimes {\bf 1},\quad \sigma=
\sigma_x\otimes {\bf 1}.
\end{equation}

\subsection{The Adiabatic Theorem}

\noindent {\bf {Theorem\  IV.1}}
 Let $H_D(s;m,d,f,\alpha)=V(s) H_D(m,d,f,\alpha )V^\dagger(s), \
s\in[0,1]$ be the family of time dependent  Dicke models with $f$ square
integrable, with square root singularity at $k=0$;
$m>\alpha^2\,\langle f|\frac{1}{|k|}|f\rangle$;  $d\ge 3$ and
$V(s)=\exp \, (i\,s\,\sigma)$ as in Eq.~(\ref{v}).  Then, $U_A$, the
adiabatic evolution associated with the ground state of
$H_D(s;m,d,f,\alpha )$, and
$U_\tau$, the Schr\"odinger evolution are close in the sense that
\begin{equation}
\Vert (U_A(s)- U_\tau(s))P(0)\vert \le
C\left\{\begin{array}{lr}\frac{1}{\tau}&\mbox{if $d>3$}\\
\frac{\sqrt{\log \tau}}{\tau}& d=3.\end{array}\right.
\end{equation} The time scale is determined by $m-\alpha^2 {\cal E}$ and
coincides with the gap without photons, $m$, up to a correction by the Lamb
shift, $\alpha^2 {\cal E}$.

 Proof: From Corollary 2.1  we find $K=\sigma\otimes\Omega$. We
will first show that a solution of the commutator equation,
Eq.~(\ref{commute}), for $d>0$, is
\begin{equation} X=\frac{i\,X_1-X_2(g)}{m-\alpha^2\,{\cal E}},\quad
Y=0,\label{x}
\end{equation} where
\begin{equation} {X}_1 =J\otimes \Omega,
\quad X_2(g) =  P\otimes(
  a^{\dagger}(g)\,\Omega + {\rm h.c}), \quad g= i\,\alpha^{\frac{3}{2}}\,
\frac{f}{|k|}.\label{g}
\end{equation} Note that the gap of the two level system $m$ is
renormalized to
${m+i\,\alpha\,\langle f|g\rangle}$, which is just the Lamb shift (See
appendix). This is a small correction, of order $\alpha^2$.

A useful formula we shall need is
\begin{equation} E\, a^\dagger(g) \Omega =a^\dagger(|k|g)\Omega.
\end{equation}

 Let us compute the commutators of ${X}_1$, ${X}_2$ with $H$:
\begin{eqnarray} [H,{X_1}] &=&
\left[\left( \begin{array}{cc} \frac{E}{\alpha}&
\sqrt{\alpha}\,{a^\dagger(f)}
 \\ \sqrt{\alpha}\, {a(f)}&m + \frac{E}{\alpha} \end{array} \right),
\left( \begin{array}{cc} {0}&-\Omega \\ \Omega & {0} \end{array}
\right)\right] \nonumber \\ &=&m\,\sigma\otimes \Omega+\sqrt{\alpha}\,
P\otimes
\Big(a^\dagger(f) \Omega+\Omega a(f)\Big).
\end{eqnarray} For the second commutator
\begin{eqnarray} [H,{X_2}] &=& \left[\left( \begin{array}{cc}
\frac{E}{\alpha}& {\sqrt{\alpha}\,a^\dagger(f)}
 \\ \sqrt{\alpha}\, {a(f)}& m + \frac{E}{\alpha} \end{array}
\right),\left ( \begin{array}{cc} { a^{\dagger}(g)\,\Omega + {\rm h.c}}&
{0} \\ {0}& 0 \end{array} \right)\right] \nonumber \\
&=&\frac{1}{\alpha}\, P\otimes  (a^{\dagger}(|k| g)\Omega - \Omega a(|k|
g)+\sqrt{\alpha}\,\left(\begin{array}{lr}0&-\langle g|f
\rangle\\ \langle f|g \rangle&0\end{array}\right)\otimes
\Omega.\label{second}
\end{eqnarray} So, if we take $g$ of Eq.~(\ref{g}) then
\begin{equation} [H,i\,{X_1}-{X_2}]=i\,\left(m-\alpha^2\,{\cal E}\right)
\sigma\otimes\Omega.
\end{equation} We see that  we can formally solve for the commutator
equation, Eq.~(\ref{commute}) provided ${\cal E}$ is finite.

This is, however, not the only condition. $X$ is a bounded operator in the
Hilbert space provided $g\in L^2$, for otherwise $a^\dagger(g)$ is ill
defined:
\begin{equation}
\int{\frac{|f|^2}{|k|^2}d^dk}\sim \int{\frac{1}{|k|^3}d^dk}<\infty.
\end{equation} The integral is finite if $d\ge 4$ but is logarithmically
divergent in $d=3$. For $d=3$ we need to squeeze $X$ back to the bounded
operators. We do that by allowing for $Y\neq 0$.

 Let $\chi_\varepsilon$ be the characteristic function of a ball of radius
$\varepsilon$ and $\chi_\varepsilon^c=1-\chi_\varepsilon$ and let
$g_\varepsilon^c =\chi_\varepsilon^c g$ and  $g_\varepsilon
=\chi_\varepsilon g$. Let us take $X_2(g_\varepsilon^c)$, which is well defined
and its norm is
$O(\alpha^{\frac{3}{2}}\sqrt{\vert\log\varepsilon\vert} )$. For $X$ we take, as before,
\begin{equation}
X=\frac{i\,X_1-X_2(g^c_\varepsilon)}{m+i\sqrt{\alpha}\langle
f|g^c_\varepsilon\rangle }.
\end{equation} From this
\begin{equation}
\Vert X\Vert = O\left(\frac{1
+\alpha^{\frac{3}{2}}|\log\varepsilon|^{1/2}}{|m-\alpha^2 {\cal E} |}.\right)
\end{equation}
 For $Y$ we take
\begin{eqnarray}
 (m+i\, \sqrt{\alpha} \langle f| g^c_\varepsilon
\rangle)\,Y&=&[H,X_2(g)-X_2(g_\varepsilon^c)]=[H,X_2(g_\varepsilon)]
\\ &=&\frac{1}{\alpha}\, P\otimes  \Big(a^{\dagger}(|k| g_e)\Omega -
 \Omega a(|k| g_e)\Big)+\sqrt{\alpha}\,\left(\begin{array}{lr}0&-\langle
g_e|f
\rangle\\ \langle f|g_e \rangle&0\end{array}\right)\otimes\Omega,\nonumber
\end{eqnarray} and we used the computation of the commutator
Eq.~(\ref{second}). With $f$ having a square root singularity,
\begin{equation}
\Vert Y\Vert = O\left(\frac{\sqrt{\alpha}
\,\varepsilon +\alpha^2\,\varepsilon}{|m-\alpha^2{\cal E} |}\right) .
\end{equation}
This puts us in the frame of theorem II.1,
except for the minor modification the log appears with a square root.
Chasing the square root establishes the main
result.\hfill
$\hfill\vbox{\hrule height 0.6pt\hbox{\vrule width 0.6pt height 1.8ex \kern
1.8ex \vrule width 0.6pt} \hrule height 0.6pt}$

 \section*{Acknowledgments} We are grateful to I.M. Sigal, S. Graffi, A.
Ori and C. Brif for useful discussions. This work was partially supported
by a grant from the Israel Academy of Sciences, the Deutsche
Forschungsgemeinschaft, and by the Fund for Promotion of Research at the
Technion.

\appendix
\section{Resonance and Lamb Shift of the Dicke Model}

The ${\cal N}=1$ sector of the Dicke model has a resonance that serves to
define the Lamb shift. The resonance is a solutions of the analytically
extended  eigenvalue equation, Eq.~(\ref{res}), see
\cite{friedrichs,howland,ms}. The real part of the shift is, by
definition, the Lamb shift of the model, while the imaginary shift is the
life time. For $d\ge 3$, the Lamb shift is dominant and the life time is a
higher order in $\alpha$. For the application to the adiabatic theorem we
need only the dominant contribution, i.e. only the Lamb shift. Computing
the Lamb shift is easy. Computing the life time is harder.  For the sake
of completeness we compute both, even though we only need one.

The eigenvalue equation is
\begin{equation}
 E-m=\alpha^2\, G(\alpha E)\label{res},
\end{equation} where $G(e)$ is defined as the analytic continuation from
the upper half plane of
\begin{equation}
\  G(e)=\int_{{\rm I\kern-.2em R}^d}\,\frac{|f|^2}{e-{|k|}}
\, d^d k,\quad \Im\, e\ge 0.
\end{equation} By taking the imaginary part, it is easy to see that
Eq~(\ref{res}) has no solution in the upper half plane. To solve the
equation in the lower half plane one needs an explicit expression, at
least for small $\alpha$, and $e$ near $\alpha m$, of this analytic
continuation. Then, we can solve Eq.~(\ref{res}) by iteration, and to
lowest order we have
\begin{equation}
 E_r\approx m+\alpha^2\, G(\alpha m).\label{iterate}
\end{equation} Clearly $G(\alpha m)\to - {\cal E}$, in the limit
$\alpha\to 0$,  so to leading order
\begin{equation}
 E_r\approx m-\alpha^2\, {\cal E}.
\end{equation} To this order, one does not see the imaginary part of the
resonance energy.
$\alpha^2\,{\cal E}$ is, by definition, \cite{dicke}, the Lamb shift of the
model. It may be worthwhile to point out that the Lamb shift for the
Hydrogen atom, \cite{bethe}, is actually of {\em higher} order, namely,
$\alpha^3
\log(\alpha^{-1})$.  Since the Lamb shift of Hydrogen also involves an
ultraviolet regularization, while the present model is ultraviolet
regular, it is not surprising that the order of the two is different. What
is surprising is that the order of Hydrogen is higher rather than lower.

Estimating the life time is, as we noted, irrelevant to the adiabatic
theorem. So a reader will loose little by skipping the rest of this
Appendix. However, for the benefit of the reader who is interested in how
the computation of the life time goes, it is given below.

 We shall show below that for $d\ge 3$, and  $|e-\alpha m|<\alpha m$, the
analytic continuation of $G(e)$ to the lower half plane, and to the next
relevant order in $\alpha$, is given by
\begin{equation}
 G(e)=-{\cal E}-i\pi \,K\Omega^d e^{d-2}\,,\quad \Im\, e\le
0,\label{analytic}
\end{equation} where $K$ is as in Eq.~(\ref{K}), and $\Omega^d$ is the
surface area of the unit ball in d dimensions. From Eq.~(\ref{iterate}),
and taking into account Eq.~(\ref{K}), we get for the Lamb shift and the
life-time:
\begin{eqnarray} E_r&\approx& m -\alpha^2 {\cal E} -i
\alpha^2 \pi \,K\,\Omega^d (m\alpha)^{d-2} \nonumber \\ &=& m -\alpha^2
{\cal E} -i
\alpha^2 \pi \Omega^d (m\alpha)^{d-1} \left\vert f(\alpha m)\right\vert^2.
\end{eqnarray}
 The  life time is  higher  order in $\alpha$ than the Lamb shift,  and is
of order $\alpha^d$. For $d=3$ this is, indeed, the order of the life time
of atomic levels that decay by dipole transition.  For small $\alpha$  the
Lamb shift dominates the life time, both in the Dicke model and in
Hydrogen.

It remains to show that the analytic continuation of $G(e)$ to the lower
half plane in a neighborhood of $m\alpha$, is indeed given by
Eq.~(\ref{analytic}). This can be done as follows: Let $B_r$ be a ball of
radius $r=2 m\alpha$. Then, in the upper half plane
\begin{equation}
\  G(e)=\left(\int_{B_r}+\int_{B_r^c}\right) \frac{|f|^2}{e-{|k|}}
\, d^d k=G_r(e)+G^c_r(e).
\end{equation} Clearly, $G^c_r(e)$ extends analytically to a half circle
in the lower half plane $|e-\alpha m|<\alpha m$. In the limit of $\alpha
\to 0$, by continuity,
\begin{equation} G^c_r(0)\to -{\cal E}.
\end{equation} This is the dominant piece, and it is real.

Consider the analytic continuation of $G_r(e)$ for $|e-\alpha m|\le \alpha
m$. Since, for small argument $f(k)$ is given by Eq.~(\ref{K}), one has
(in the upper half plane)
\begin{equation} G_r(e)=K\Omega^d\int_0^{2m\alpha} \frac{k^{d-2}}{e-k}\,
dk= K\Omega^d\int_\gamma \frac{k^{d-2}}{e-k}\, dk,
\end{equation} where $\gamma$ is the obvious semi-circle in the complex
$k$ plane and $\Omega^d$ the surface area of the unit ball in $d$
dimensions. The right hand side is analytic in $e$ in the lower half plane
provided $|e-\alpha m|<\alpha m$, and so gives the requisite analytic
continuation. Since $e$ is small, and of order
$\alpha$, to leading order, we have
\begin{eqnarray} G_r(e)&=&K\Omega^d\int_\gamma
\frac{(k-e +e)^{d-2}}{e-k}\,dk\nonumber \\ &=&-K\Omega^d\,
\sum_{j=0}^{d-2} \left(\begin{array}{c} d-2\\
j\end{array}\right)\,\,e^{d-j-2}\int_\gamma (k-e)^{j-1}\,dk
\nonumber \\ &\approx&-K\Omega^d\, \,e^{d-2}\int_\gamma
\frac{dk}{k-e}\nonumber \\ &=&-\, K\Omega^d\, e^{d-2}
 \left(i\,\pi\,+\log\left(\frac{2\alpha
m}{e}\right)\,+\,O(\alpha\log\alpha)\right)
\end{eqnarray} and the error term in approximation that we did not compute
is real and being sub-dominant to ${\cal E}$ is irrelevant.

\end{document}